\documentclass{svmult}

\usepackage{makeidx}         
\usepackage{graphicx}        
\usepackage{multicol}        
\usepackage[bottom]{footmisc}
\usepackage{amsmath,amssymb}

\makeindex             

\newcommand{\be}{\begin{equation}}
\newcommand{\ee}{\end{equation}}
\newcommand{\ba}{\begin{array}}
\newcommand{\ea}{\end{array}}

\newcommand{\pref}[1]{(\ref{#1})}

\def\fn1{f^{n+1}}

\begin{document}

\title*{Increasing efficiency through optimal RK time integration of
  diffusion equations}
\titlerunning{Optimal RK time integration of
  diffusion equations}
\author{Fausto Cavalli\inst{1} \and Giovanni Naldi\inst{1}  \and Gabriella Puppo\inst{2}  \and Matteo
  Semplice\inst{1}}
\authorrunning{F. Cavalli, G. Naldi, G. Puppo, M. Semplice}
\institute{Dipartimento di Matematica, Universit\`a di Milano, via
  Saldini 50, 20133 Milano, ITALY.
  \texttt{\{cavalli,naldi,semplice\}@mat.unimi.it} 
  \and 
  Dipartimento di Matematica, corso Duca degli Abruzzi, 24, 10129 Torino, Italy
  \texttt{gabriella.puppo@polito.it}
}
%
%

\maketitle

\begin{abstract}
\noindent The application of Runge-Kutta schemes designed to enjoy a
large region of absolute stability can significantly increase the
efficiency of numerical methods for PDEs based on a method of lines
approach. In this work we investigate the improvement in the
efficiency of the time integration of relaxation schemes for
degenerate diffusion problems, using SSP Runge-Kutta schemes and
computing the maximal CFL coefficients. This technique can be extended
to other PDEs, linear and nonlinear, provided the space operator has
eigenvalues with a non-zero real part.
\end{abstract}

\section{Introduction}\label{sect:introduction}
The integration of evolution PDE's through a method of lines approach
leads to the solution of large systems of ODE's. Often such ODE's are
stiff or moderately stiff; therefore the possibility of increasing the
stability region of the time integrator can lead to a significant
increase in efficiency for explicit or semi-implicit time-integration.

Specifically, we consider the system of PDE's 
\be 
u_t + f_x(u) = D p_{xx}(u), \label{eq:general_system}
\ee
where $f(u)$ is hyperbolic, i.e. the Jacobian of $f$ is provided with
real eigenvalues and a basis of real eigenvectors for each $u$, while
$p(u)$ is a non decreasing 
Lipshitz continuous function, with Lipshitz constant $\mu$. We assume
the system has been fully 
discretized in space on a grid of $N$ points $x_j$, $j=1,\dots N$, and
we denote with $U(t)=[U_1(t), \dots, U_N(t)]^T$ the vector of the grid
values of the numerical solution at time $t$.
The space discretized system can be written in the form: 
\be
\frac{dU}{dt} = L(U(t)),  \label{eq:general_semidiscrete}
\ee
leading to an autonomous system of $N$ non linear first order
ODE's. If we consider for instance a system of conservation laws,
$D=0$, the operator $L$ obtained from a conservative space
discretization, will be written as: 
\[
L(U) = -\frac1h \left( F_{j+1/2} - F_{j-1/2} \right),
\]
where $F_{j+1/2}$ is the numerical flux consistent with the physical
flux $f(u)$ in the usual sense of the Lax-Wendroff theorem and $h$ is
the grid spacing. 

In recent years, much research has focused on the efficient
integration of the semidiscrete system
\pref{eq:general_semidiscrete}. In particular, in this work we will
concentrate on the performance of optimal Runge-Kutta schemes,
characterized by a large stability region, introduced in
\cite{SR02}. 

Optimal Runge-Kutta schemes are built choosing an accuracy order $p$
and a number of stages $s$, with $s \ge p$. In principle, once $s$ and
$p$ are fixed, the coefficients of the Butcher tableaux defining the
Runge-Kutta schemes are computed maximizing in some sense the
stability region and keeping as constraints the fulfilment of the
accuracy requirements. To compute their optimal schemes, Spiteri and
Ruuth in \cite{SR02} start from a strong stability condition which
requires that the operator $L$ be non linearly stable with respect to
a suitable norm for a certain CFL with Forward Euler integration,
namely:
\[
|| U^n + \Delta t L(U^n)|| \leq ||U^n||, \quad \forall \Delta t \leq \Delta t_{FE}.
\]
Once this assumption is satisfied, the optimal schemes proposed in
\cite{SR02} do yield considerable savings in CPU time for a given
accuracy.

The idea is that an $s$ stages Runge-Kutta scheme applied to
\pref{eq:general_semidiscrete} can be written as a convex combination
of $s$ Forward Euler steps as (see also \cite{GST01}): 

\begin{eqnarray}
U^{(1)} & = & U^n  \\
U^{(i)} & = & \sum_{k=1}^{i-1} \alpha_{ik} \left[ U^{(k)} + \Delta t
  \frac{\beta_{ik}}{\alpha_{ik}} L(U^{(k)}) \right] \\ 
U^{(n+1)} & = & U^{(s)}.  \label{eq:general_RK}
\end{eqnarray}
Thus if the space discretized operator $L$ is strongly stable for
$\Delta t \leq \Delta t_{FE}$ with Forward Euler time integration,
than the scheme in \pref{eq:general_RK} will be strongly stable for: 
\be
\Delta t \leq \lambda \, \Delta t_{FE}, \qquad \lambda = \min \frac{\alpha_{ik}}{\beta_{ik}}
                  \label{eq:cfl_fe}
\ee
Note that \pref{eq:cfl_fe} is only a sufficient condition for
stability: clearly it may be possible to violate \pref{eq:cfl_fe} and
still find a stable scheme. 

The problem is that several high order space discretization operators
$L$ {\em are not} stable under the Forward Euler scheme, and therefore
it is not possible to use the estimate \pref{eq:cfl_fe} to guarantee
that optimal SSP schemes will  improve the efficiency of the time
integration of \pref{eq:general_semidiscrete}. 

In the case of pure convection, that is for $D=0$ in
\pref{eq:general_system}, high order space discretizations of $f_x(u)$
have purely imaginary eigenvalues $\nu$ up to high order powers of the
mesh width $h$, that is $\Re(\nu) = O(h)^p$ and therefore are not
stable under Forward Euler time integration. This is the case for
instance of the fifth order WENO space discretization, see also
\cite{Spiteri:private}, but we conjecture that the same holds for
other widely used high order space discretizations for convective
operators, such as ENO.

On the other hand, when $D\neq0$, the eigenvalues of the semidiscrete
operator $L$ do have a negative real part of order $1$, and therefore
can be made stable under Forward Euler for a non zero $\Delta t_{FE}$.
In this work we study the stability of a family of high order
numerical fluxes coupled with optimal RK-SSP schemes for
\pref{eq:general_system}, in the case of pure diffusion, i.e. $f
\equiv 0$. In this case we find that several widespread space
discretization schemes are stable with Forward Euler time integration,
and therefore optimal RK-SSP schemes do yield a significant increase
in the allowable CFL. Since in this case the eigenvalues of the exact
operator are real, we even find that the stability estimate in
\pref{eq:cfl_fe} may be quite pessimistic, because it underestimates
the stability region of some SSP schemes. In these cases, it is easy
to compute numerically the maximal CFL, obtaining a further
improvement in the efficiency of the scheme.

\section{Diffusive relaxation}\label{sect:diffusive}

We consider high order approximations of the degenerate parabolic equation
\be 
u_t  = D p_{xx}(u), \label{eq:parabolic_system}
\ee
using relaxation schemes, a technique initially proposed in \cite{JX95}. Following
\cite{NP00} and \cite{arxiv0604572}, we rewrite the system as a hyperbolic
system with a stiff source term depending on a parameter
$\varepsilon$, which formally relaxes on the original parabolic equation as
$\varepsilon \rightarrow 0$, namely: 
\begin{equation} \label{3eq}
\left\{
\begin{array}{ll}
\displaystyle
\frac{\partial u}{\partial t} + \frac{\partial v}{\partial x} = 0 \\
\\
\displaystyle \frac{\partial v}{\partial t} + \Phi^2 \frac{\partial
  w}{\partial x}  = 
    -\frac1\varepsilon v
    +\left(\Phi^2- \frac{D}{\varepsilon}\right) \frac{\partial
  w}{\partial x} \\
\\
\displaystyle
\frac{\partial w}{\partial t} + \frac{\partial v}{\partial x} =
    -\frac1\varepsilon (w-p(u))
\end{array} 
\right.
\end{equation}
where $\Phi^2$ is a suitable positive parameter. Formally, as
$\varepsilon \rightarrow 0$, $w \rightarrow p(u)$, and $v \rightarrow
- \partial_x p(u)$ and the original equation \pref{eq:parabolic_system} is
recovered. We integrate the system above with an IMEX Runge-Kutta
scheme \cite{PR05}. In this fashion, the stiff source term is
implicit and does not require restrictive stability conditions, while
the linear convective term is explicit.  

In this work we consider only the relaxed scheme, which is obtained
setting $\varepsilon = 0$ in the discretized equations. Let $\tilde{a}_{i,k}$ and $\tilde{b}_i$ be the coefficients forming the Butcher tableaux of the explicit Runge-Kutta scheme in the IMEX pair. The
computation of the first stage value of the Runge Kutta scheme reduces
to:  
\be
\ba{l} u^{(1)}=u^n \\ w^{(1)}=p(u^n) \\ v^{(1)}=-D\partial_x w^{(1)}\ea.
\ee
For the following stages the first equation is
\be\label{originalRK}
 u^{(i)}=u^n -\Delta t \sum_{k=1}^{i-1}\tilde{a}_{i,k}
   \partial_x v^{(k)}.
\ee
In the other equations the convective terms are dominated by the
source terms and thus $v^{(i)}$ and $w^{(i)}$ are given by  
\be
 w^{(i)}=p(u^{(i)}), \qquad v^{(i)}=-D\partial_x w^{(i)}.
\ee
In this fashion, due to the particular structure of the relaxation
scheme we are considering, only the explicit tableaux of the
Implicit-Explicit Runge-Kutta scheme enters the actual computation,
while the coefficients of the implicit scheme drop out as the
relaxation step is computed. For this reason we can apply any explicit
Runge-Kutta scheme to advance in time the solution $u$. Finally the updated solution is given by: 
\be
u^{n+1} = u^n  - \Delta t \sum_{i=1}^{\nu} \tilde{b}_i \partial_x v^{(i)}
\ee 
The overall accuracy of the scheme depends on the accuracy of the
numerical flux used to approximate $\partial_x v^{(i)}$ and on the
accuracy of the Runge-Kutta explicit scheme.  

A non linear stability analysis for the first order version of this
scheme (upwind numerical flux on a piecewise constant space reconstruction to evaluate $\partial_x v$ and Forward Euler in time) yields a parabolic stability restriction of the form: 
\be
\Delta t \leq \frac{2h^2}{\mu} \frac{1}{1+2h \phi},
\ee
where $\mu$ is the Lipshitz continuity constant of $p(u)$, see
\cite{arxiv0604572}. The parabolic CFL implies that to  obtain a scheme of order $p$ one should use a $p$th order
accurate numerical flux for $\partial_x v$ and a $p/2$ order accurate Runge-Kutta method in time.

A linear stability analysis of several higher order numerical fluxes
yields again a parabolic CFL, of the form $ \Delta t \leq C_1
h^2(1-C_2h\Phi)/\mu$,  with
constant $C_1$ given by Table \ref{table_CFL}.
Note that the scheme is unstable for $\Delta t=C_1h^2/\mu$, but for a
sufficiently small $h$ it is enough to pick $\Delta t=(C_1-\delta)h^2/\mu$
for a suitable positive $\delta$, i.e. the scheme is stable provided
$C_1$ is slightly decreased, see Table \ref{tab:matteo}. Note that the
stability requirement becomes more strict as space accuracy increases,
while it loosens as time accuracy increases. 

\begin{table}
\begin{center}
\begin{tabular}{|c|c|c|c|}
\hline
                & RK1 & RK2 & RK3 \\
\hline
P-wise constant & 2   & 2   & 2.51\\
\hline
P-wise linear   & 0.94&0.94 & 1.18\\
\hline
WENO5           & 0.79&0.79 & 1   \\
\hline
\end{tabular}
\end{center}
\caption{CFL constant $C_1$ for a few space reconstruction algorithms (from linear analysis) for standard first, second and third order RK schemes}  \label{table_CFL}
\end{table}

\section{Diffusive relaxation and SSP schemes}

We start fixing a standard notation for $s$ stages Runge-Kutta schemes. As in
\cite{SR02} we denote SSP(s,p) the optimal strongly stable
Runge-Kutta scheme of order $p$ with $s$ stages. We point out that
SSP(1,1) is the Forward Euler scheme, SSP(2,2) the Heun scheme
and SSP(3,3) the TVD third order Runge-Kutta method of \cite{SO88},
which is probably the high order Runge-Kutta scheme most frequently 
used in conjunction with high order space discretizations for
\pref{eq:general_system}.

\begin{table}
\begin{center}
\begin{tabular}{|c|c|c|c|c|c|} 
\hline
     & s=1 & s=2 & s=3 & s=4 & s=5\\ \hline
 p=1 & 1   & 2   & 3   & 4   & 5  \\   \hline
 p=2 &     & 1   & 2   & 3   & 4  \\   \hline
 p=3 &     &     & 1   & 2   & 2.65  \\   \hline
\end{tabular}
\caption{Improved stability coefficients for several $s$ stages order $p$
  Runge-Kutta schemes (from \cite{SR02})}  
\label{table:optimal_cfl}
\end{center}
\end{table}

From the first column of Table \ref{table_CFL} we note that all the
numerical fluxes considered are at least linearly stable with the
Forward Euler scheme. Thus the theory of strongly stable Runge-Kutta
schemes can be applied in this case. In particular we consider the
schemes introduced in \cite{SR02}, for which the improved stability
coefficients $\lambda$ of \pref{eq:cfl_fe} can be found in Table
\ref{table:optimal_cfl}. To indicate the fact that these coefficients
are found through the theory of strongly stable Runge Kutta schemes
and to identify to which scheme they apply, we will label these
coefficients as $\lambda_{\mbox{\tiny SSP}}(s,p)$.

\begin{figure}
\begin{center}
\includegraphics[width=.4\textwidth]{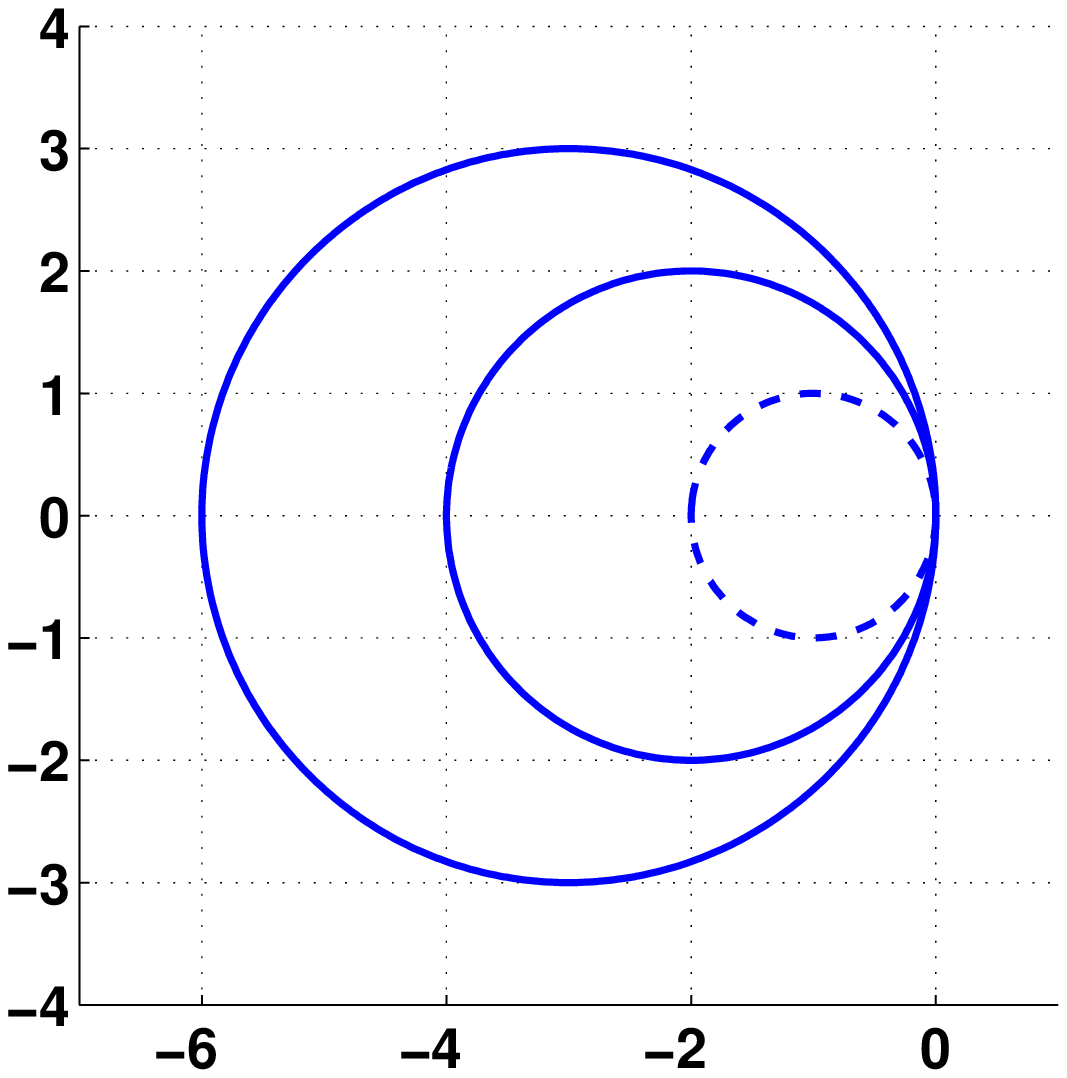}
\hfil
\includegraphics[width=.4\textwidth]{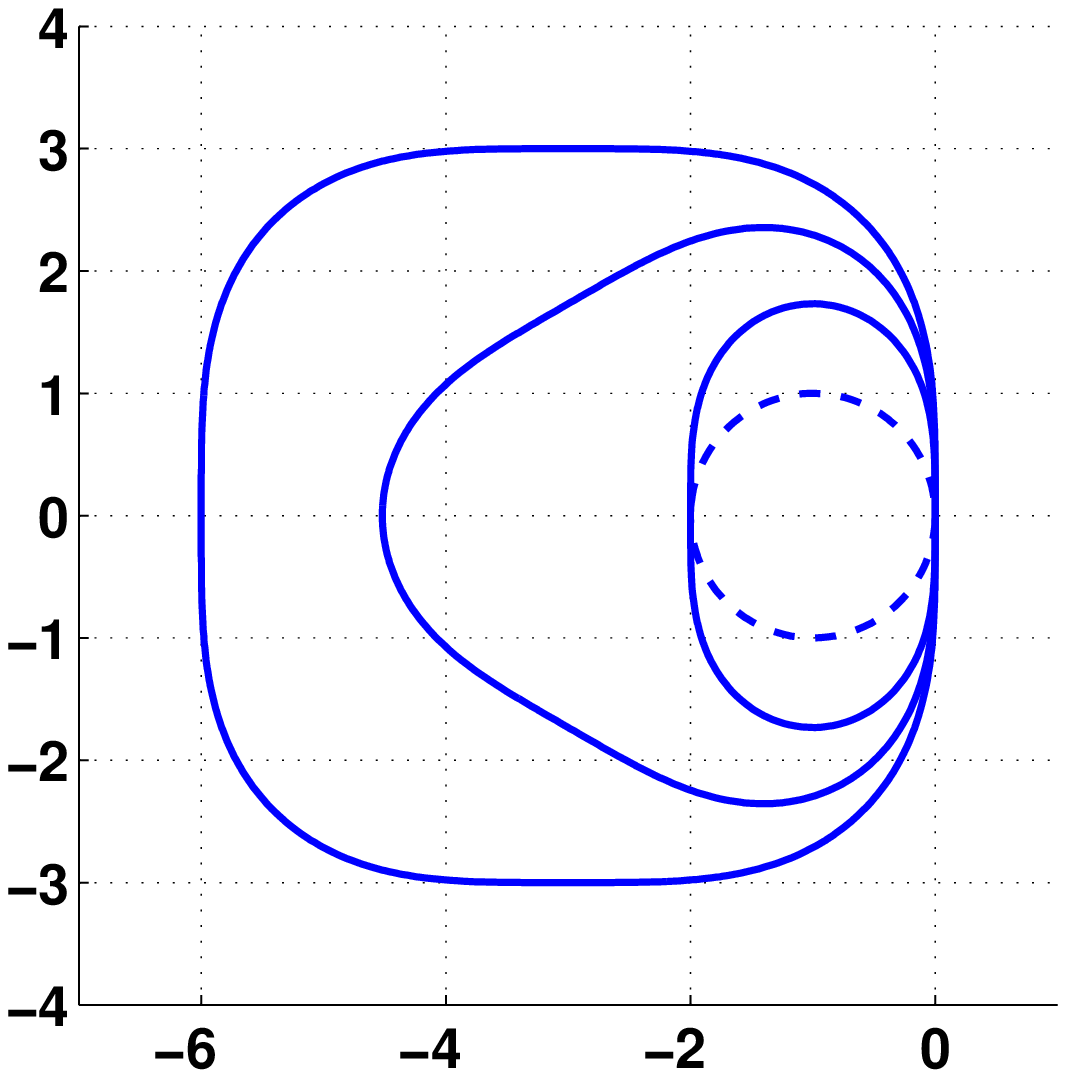}\\
\includegraphics[width=.4\textwidth]{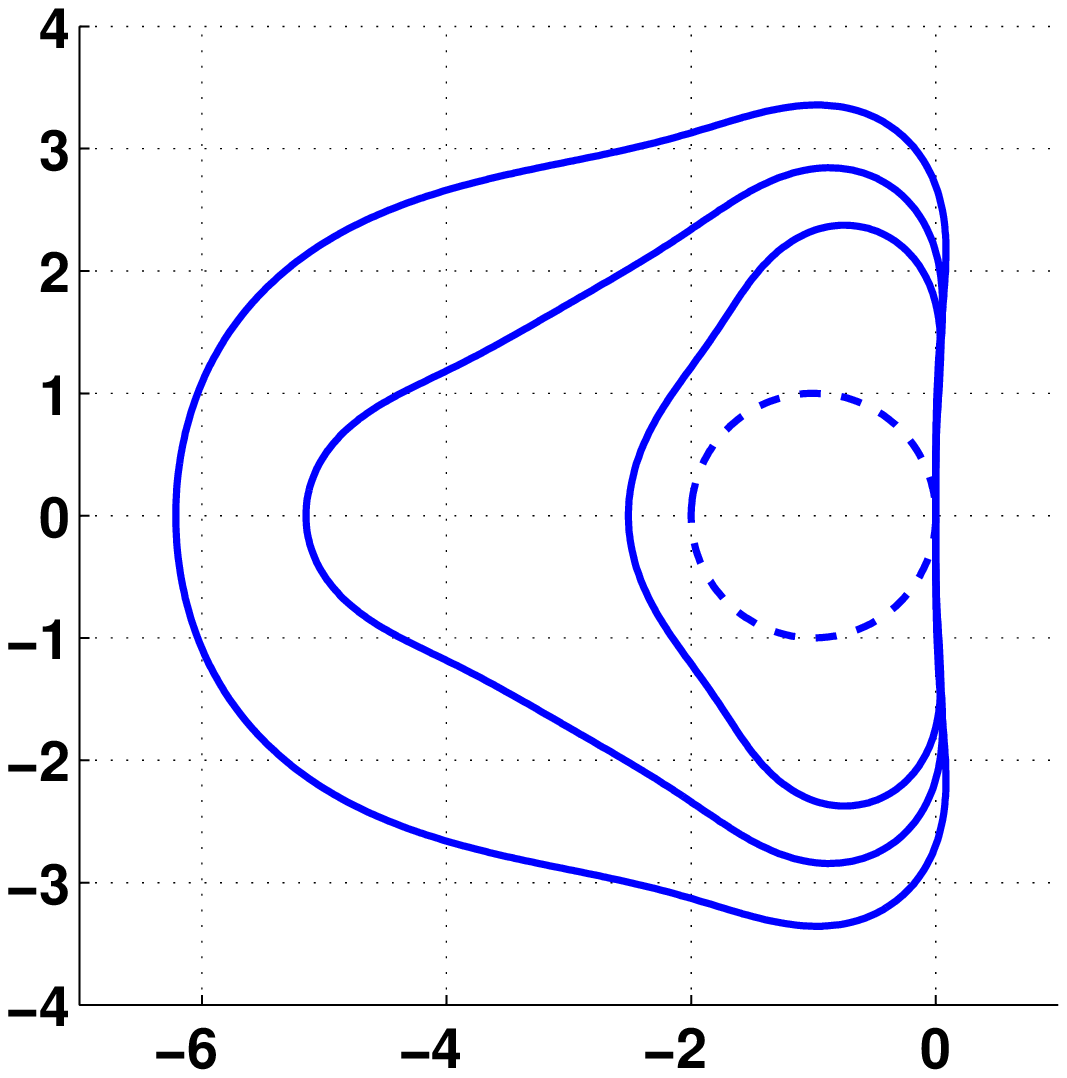}
\end{center}

\caption{Regions of absolute stability for optimal SSP schemes. On the
  top left: order 1, stages 1,2,3. On the top right: order 2, stages
  2,3,4. On the bottom: order 3, stages 3,4,5. The dashed line is the
  Forward Euler scheme.}
\label{fig:fausto}
\end{figure}

The coefficients $\lambda_{\mbox{\tiny SSP}}(s,p)$ are not optimal in
the case of pure diffusion. They can be further improved considering
the actual stability region of SSP schemes and recalling that for the
purely diffusive operator we are considering, only the intersection of
the stability region with the real line is relevant.

In Figure \ref{fig:fausto} we find the regions of absolute
stability for some SSP schemes 
given in \cite{SR02}. On the top left of the figure we find the stability
regions of schemes of
order $1$ with $1,2$ and $3$ stages. We point out that here the improvement
in the CFL constant is exactly balanced by the increased computational
effort due to the higher number of stages: there is no gain in
efficiency with respect to the standard Forward Euler scheme.
In the top right of Figure \ref{fig:fausto} there are the stability
regions of schemes of
order $2$ with $2,3$ and $4$ stages and, for comparison, of the
Forward Euler scheme (dashed line). The SSP theory underestimates the effective CFL
and direct inspection of the stability plot suggests that the
stability coefficient $\lambda$ appearing in Table \ref{tab:matteo} can be
increased finding the intersection of the stability curve with the
real axis.
Finally, on the bottom of Figure \ref{fig:fausto} we show the stability
regions of schemes of
order $3$ with $3,4$ and $5$ stages. Again the CFL gain of
Table \ref{table:optimal_cfl} can be improved by direct inspection of
the graph. 

Let $\eta_{s,p}$ be the abscissa of the
intersection of the stability curve with the negative real axis. Then the
maximal gain in CFL with respect to Forward Euler, for a problem with
real eigenvalues, is given by
\[ \lambda_{\mbox{\tiny OPT}}(s,p) = \frac{|\eta_{s,p}|}{|\eta_{1,1}|} 
         = \frac{|\eta_{s,p}|}{2} 
\]
For several schemes in the figure it is easy to see that
$\lambda_{\mbox{\tiny OPT}}(s,p)>\lambda_{\mbox{\tiny SSP}}(s,p)$.

The optimal CFL number for a given scheme can be found multiplying the
coefficient $C_1-\delta$ determined by the space discretization by
the proper stability coefficient $\lambda$, i.e.
\(
\Delta t \leq (C_1-\delta) \lambda_{\mbox{\tiny OPT}}(s,p) h^2/\mu
\).

\begin{table}
\begin{center}
\begin{tabular}[b]{|c|c|c|c|}
\multicolumn{4}{l}{SSP(2,s) + WENO5}\\ \hline
stages & CFL & order & $N_f$ 
\\ \hline
2 & $0.78$        & $4$ & $810$ 
\\ \hline
3 & $2\times0.78$ & $4$ & $606$ 
\\ \hline
4 & $3\times0.78$ & $4$ & $540$ 
\\ \hline \hline
\multicolumn{4}{l}{SSP(3,s) + WENO5}\\ \hline
stages & CFL & order & $N_f$ 
\\ \hline
3 & $0.78$        & $4.5$ & $1230$ 
\\ \hline
4 & $2\times0.78$        & $5.2$ & $820$ 
\\ \hline
5 & $2.65\times0.78$        & $5.2$ & $770$ 
\\ \hline
\end{tabular}
\hfil
\begin{tabular}[b]{|c|c|c|c|}
\multicolumn{4}{l}{SSP(2,s) + WENO5}\\ \hline
stages & CFL & order & $N_f$ 
\\ \hline
2 & $1\times0.78$ & $4$ & $810$ 
\\ \hline
3 & $\mathbf{2.259\times}0.78$ & $4$ & $537$ 
\\ \hline
4 & $3\times0.78$ & $4$ & $540$ 
\\ \hline \hline
\multicolumn{4}{l}{SSP(3,s) + WENO5}\\ \hline
stages & CFL & order & $N_f$ 
\\ \hline
3 & $\mathbf{1.256\times}0.78$   & $4.8$ & $978$ 
\\ \hline
4 & $\mathbf{2.574\times}0.78$   & $5.6$ & $636$ 
\\ \hline
5 & $\mathbf{3.106\times}0.78$   & $5.4$ & $660$ 
\\ \hline
\end{tabular}
\end{center}
\caption{Order of convergence and number of numerical flux function
  evaluations $N_f$ (with 80 grid points) for the SSP-CFL      $\lambda_{\mbox{\tiny{SSP}}}(s,p)$ (left) and the
  maximal CFL $\lambda_{\mbox{\tiny{OPT}}}(s,p)$ (right) for several RK schemes.}
\label{tab:matteo}
\end{table}

To measure the computational complexity of a scheme we compute the number $N_f$
of numerical flux evaluations needed to reach a fixed integration
time. Thus, for a given numerical flux, the most efficient scheme has
the lowest value of $N_f$. 
Table \ref{tab:matteo} shows the different values of $N_f$ obtained
with CFL chosen according to the values of 
$\lambda_{\mbox{\tiny SSP}}(s,p)$ on the left and 
$\lambda_{\mbox{\tiny OPT}}(s,p)$ on the right. 
The grid spacing $h$
is the same for all values of $N_f$, namely, $h=1/80$. The table
contains also data on the accuracy of the space-time scheme. The
accuracy was evaluated using four different grids, and modifying the
time step according to $h$ and the chosen value of $\lambda$.   
We don't show data for first order schemes since 
$\lambda_{\mbox{\tiny SSP}}(s,1)=\lambda_{\mbox{\tiny OPT}}(s,1)$. The
table shows the increased efficiency of the SSP(s,p) schemes with
$s>p$ (left column) and the gain obtained by a better estimate of
$\lambda$ (right column). In the table, the optimal values of
$\lambda$ are indicated in bold face, when they are sensitively larger
than the corresponding $\lambda_{\mbox{\tiny SSP}}(s,p)$. In
particular, the improvement between the standard SSP(3,3) and the five
stages SSP(5,3) with the optimal CFL is quite striking. 

The data in Table \ref{tab:matteo} refer to a linear diffusion
problem, although we find analogous results for a nonlinear degenerate
diffusion equation. We performed tests on the self-similar Barenblatt
solution of the porous media equation. In this case the order of
accuracy is limited by the non regularity of the solution, but we find
that the errors with respect to the exact solution slightly decrease
using the SSP(s,p) with $s>p$ and the optimal $\lambda$.

\section{Final remarks}

We have shown that the theory of \cite{SR02} can be applied
to diffusion equations in the relaxation framework to improve the
efficiency of the time integration. Moreover the fact that the
eigenvalues of the semidiscrete operator are real numbers permits to further improve the
efficiency of the schemes. We expect that analogous results can be
obtained for convection-diffusion operators, thanks to the nonzero
real part of the eigenvalues of the discrete operators. This allows to
achieve stability under Forward Euler integration.

We also note that the same framework can be applied to other space
discretizations besides the numerical fluxes obtained via relaxation
schemes: the key factor is the localization of the eigenvalues of the
differential operator, which must have a non zero real part. On the
other hand, the plots of absolute stability regions show that
Runge-Kutta schemes with number of stages $s>p$ can be built with
improved stability regions, notwithstanding stability under Forward
Euler. These schemes can be applied to semidiscrete operators even in
the convective regime, but their stability conditions cannot be
derived from the behaviour of the operator under the Forward Euler
scheme.  

\bibliographystyle{alpha}
\bibliography{ssp}
\end{document}